# Quaternionic Salkowski Curves and Quaternionic Similar Curves

## Mehmet Önder
*Delibekirli Village, Tepe Street, No: 63, 31440, Kırıkhan, Hatay, Turkey.*
E-mail: mehmetonder197999@gmail.com

## Abstract
In this paper, we give the definitions and characterizations of quaternionic Salkowski, quaternionic anti-Salkowski and quaternionic similar curves in the Euclidean spaces $E^3$ and $E^4$. We obtain relationships between these curves and some special quaternionic curves such as quaternionic slant helices and quaternionic $B_2$-slant helices.

## 1. Introduction and Preliminaries

In the differential geometry, special curves which satisfy some relationships between their curvatures and torsions have an important role. The most popular one of these curves is general helix which is defined by the property that the tangent of the curve makes a constant angle with a fixed straight line called the axis of the general helix [1]. Moreover, recently some new special curves have been defined and studied. Izumiya and Takeuchi [2] have defined slant helix which is a special curve whose principal normal vector makes a constant angle with a fixed direction. Önder et al. [3] have considered the notion of slant helix in $E^4$ and defined $B_2$-slant helix. Furthermore, Salkowski [4] defined the curves with constant curvature but non-constant torsion by an explicit parametrization. Later, Monterde [5] has given some characterizations of Salkowski and anti-Salkowski curves. El-Sabbagh and Ali [6] have defined a new curve couple called similar curves whose arc-length parameters have a relationship and their tangents are the same. These curves have also been studied in different spaces [7,8].

In this paper, we define quaternionic and spatial quaternionic Salkowski curves, anti-Salkowski curves and similar curves. We obtain the characterizations for these special quaternionic and spatial quaternionic curves. First, we give the basic elements of the theory of quaternions and quaternionic curves. A more complete elementary treatment of quaternions and quaternionic curves can be found in references [9-11].

A real quaternion $q$ is an expression of the form $q = a_1 e_1 + a_2 e_2 + a_3 e_3 + a_4 e_4$, where $a_i$, $(1 \le i \le 4)$ are real numbers, and $e_i$, $(1 \le i \le 4)$, $e_4 = 1$ are quaternionic units which satisfy the non-commutative multiplication rules $e_i \times e_i = -e_4$, $(1 \le i \le 3)$; $e_i \times e_j = -e_j \times e_i = e_k$, $(1 \le i, j, k \le 3)$ where $(ijk)$ is an even permutation of $(123)$ in the Euclidean space. The algebra of the quaternions is denoted by $Q$ and its natural basis is given by $\{e_1, e_2, e_3, e_4\}$. A real quaternion is given by the form $q = s_q + v_q$ where $s_q = a_4$ is scalar part and $v_q = a_1 e_1 + a_2 e_2 + a_3 e_3$ is vector part of $q$.

The conjugate of $q = s_q + v_q$ is defined by $\bar{q} = s_q - v_q$. This defines the symmetric real-valued, non-degenerate, bilinear form as follows:

$$h : Q \times Q \to \mathrm{IR}, \ (q, p) \to h(q, p) = \frac{1}{2}(q \times \bar{p} + p \times \bar{q}), \qquad (1)$$

which is called the quaternion inner product. Then the norm of $q$ is

$$\|q\|^2 = h(q, q) = q \times \bar{q} = \bar{q} \times q = a_1^2 + a_2^2 + a_3^2 + a_4^2. \qquad (2)$$

If $\|q\|=1$, then $q$ is called unit quaternion. Let $q = s_q + v_q = a_1e_1 + a_2e_2 + a_3e_3 + a_4e_4$ and $p = s_p + v_p = b_1e_1 + b_2e_2 + b_3e_3 + b_4e_4$ be two quaternions in $Q$.

**Definition 1.** ([9]) $q$ is called a spatial quaternion whenever $q + \bar{q} = 0$ and called a temporal quaternion whenever $q - \bar{q} = 0$. Then a general quaternion $q$ is given as

$$q = \frac{1}{2}(q+\bar{q}) + \frac{1}{2}(q-\bar{q}).\qquad(3)$$

The quaternion $\frac{1}{2}(q-\bar{q})$ is a spatial quaternion and called spatial part of $q$ and the quaternion $\frac{1}{2}(q+\bar{q})$ is a temporal quaternion and called temporal part of $q$.

**Theorem 1.** ([9]) *The three-dimensional Euclidean space $E^3$ is identified with the space of spatial quaternions $\{q \in Q : q + \bar{q} = 0\}$ in an obvious manner. Let $I = [0,1]$ be an interval in real line $IR$ and let $\alpha : I \subset IR \to Q$, $s \to \alpha(s) = \sum_{i=1}^{3} \alpha_i(s) e_i$ be an arc-lengthed curve with nonzero curvatures $\{k, r\}$ and $\{\vec{t}(s), \vec{n}_1(s), \vec{n}_2(s)\}$ denotes the Frenet frame of the curve $\alpha(s)$. Then the Frenet formulae of the quaternionic curve $\alpha(s)$ are given by*

$$\vec{t}' = k\vec{n}_1, \quad \vec{n}_1' = -k\vec{t} + r\vec{n}_2, \quad \vec{n}_2' = -r\vec{n}_1,\qquad(4)$$

*where $k(s)$ is principal curvature and $r(s)$ is torsion of $\alpha(s)$.*

**Theorem 2.** ([9]) *The four-dimensional Euclidean space $E^4$ is identified with the space of real quaternions. Let $I = [0,1]$ be an interval in real line $IR$ and let $\gamma : I \subset IR \to Q$, $s \to \gamma(s) = \sum_{i=1}^{4} \alpha_i(s) e_i$, $e_4 = +1$ be a smooth curve in $E^4$ with nonzero curvatures $\{K, k, r-K\}$ and $\{\vec{T}(s), \vec{N}(s), \vec{B}_1(s), \vec{B}_2(s)\}$ denotes the Frenet frame of the curve $\gamma(s)$. Then the Frenet formulae of the quaternionic curve $\gamma(s)$ are given by*

$$\vec{T}' = K\vec{N}, \quad \vec{N}' = -K\vec{T} + k\vec{B}_1, \quad \vec{B}_1' = -k\vec{N} + (r-K)\vec{B}_2, \quad \vec{B}_2' = -(r-K)\vec{B}_1,\qquad(5)$$

*where $K(s)$ is principal curvature, $k(s)$ is torsion and $(r-K)(s)$ is bitorsion of $\gamma(s)$.*

$\gamma(s)$ is called a quaternionic $\vec{B}_2$-slant helix (resp. $\vec{B}_1$-slant helix) if its second principal normal $\vec{B}_2$ (resp. first principal normal $\vec{B}_1$) makes a constant angle with a fixed direction $\vec{U}$ [12].

## 2. Spatial Quaternionic Salkowski Curves in $E^3$

In this section, we introduce the explicit parametrization of a spatial quaternionic Salkowski curve in $E^3$. First, we give the following definition.

**Definition 2.** ([12]) Let $\alpha(s) : I \subset IR \to Q$ be an arc-length parametrized spatial quaternionic curve with Frenet frame $\{\vec{t}(s), \vec{n}_1(s), \vec{n}_2(s)\}$. We call $\alpha(s)$ as spatial quaternionic slant helix if $\vec{n}_1(s)$ makes a constant angle with a unit and constant real spatial quaternion.

If $\vec{n}_2(s)$ makes a constant angle with a unit and constant real spatial quaternion, then $\alpha(s)$ is called spatial quaternionic $\vec{n}_2$-slant helix.

Now, we give the definition and characterizations of spatial quaternionic Salkowski curves in $E^3$ as follows.

**Definition 3.** A spatial quaternionic curve $\alpha(t)$ with constant curvature and non-constant torsion is called spatial quaternionic Salkowski curve. For any $m \in IR$, an example of spatial quaternionic Salkowski curve $\alpha_m(t)$ is given by the parametrization

$$\alpha_m(t) = \frac{1}{\sqrt{1+m^2}}\left[\left(\frac{n-1}{4(1+2n)}\sin((1+2n)t) - \frac{1+n}{4(1-2n)}\sin((1-2n)t) - \frac{1}{2}\sin t\right)\vec{e}_1 \right.$$

$$+ \left(\frac{1-n}{4(1+2n)}\cos((1+2n)t) + \frac{1+n}{4(1-2n)}\cos((1-2n)t) + \frac{1}{2}\cos t\right)\vec{e}_2$$

$$\left. + \frac{1}{4m}\cos(2nt)\vec{e}_3 \right]$$

with $n = \frac{m}{\sqrt{1+m^2}}$. Then we have $\|\alpha'_m(t)\| = \frac{\cos(nt)}{\sqrt{1+m^2}}$, $s = \frac{\sin(nt)}{m}$, $k(t) = 1$, $r(t) = \tan(nt)$.

***Theorem 3.*** *Let $\alpha(s): I \subset IR \to Q$ be an arc-length parametrized spatial quaternionic curve with principal curvature $k \equiv 1$. The principal normal vector $\vec{n}_1(s)$ makes a constant angle $\theta$ with a fixed direction if and only if the torsion is given by $r(s) = \pm \frac{s}{\sqrt{\tan^2\theta - s^2}}$.*

**Proof:** Let $\vec{d}$ be fixed unit spatial quaternion which makes a constant angle $\theta$ with principal normal vector $\vec{n}_1(s)$. Therefore we have $h(\vec{d}, \vec{n}_1) = \cos\theta$. Differentiating that with respect to $s$ gives

$$h(\vec{d}, -\vec{t} + r\vec{n}_2) = 0. \tag{6}$$

and from Eq. (6) we have $h(\vec{d}, \vec{t}) = r h(\vec{d}, \vec{n}_2)$. If we put $h(\vec{d}, \vec{n}_2) = x$, we write $\vec{d} = rx\vec{t} + \cos\theta \vec{n}_1 + x\vec{n}_2$. Since we assume $\vec{d}$ is unit, i.e., $\|\vec{d}\| = 1$, we have

$$x = \pm \frac{\sin\theta}{\sqrt{1+r^2}}. \tag{7}$$

Then, the vector $\vec{d}$ is given by $\vec{d} = \pm \frac{r\sin\theta}{\sqrt{1+r^2}}\vec{t} + \cos\theta \vec{n}_1 \pm \frac{\sin\theta}{\sqrt{1+r^2}}\vec{n}_2$. Differentiating Eq. (6) with respect to $s$, it follows $h(\vec{d}, -(k+r^2)\vec{n}_1 + r'\vec{n}_2) = 0$. Then, from Eq. (7) we obtain the differential equation $\pm \tan\theta \frac{r'}{(1+r^2)^{3/2}} + 1 = 0$. By integration, we get

$$\pm \tan\theta \frac{r}{\sqrt{1+r^2}} + s + c = 0, \tag{8}$$

where $c$ is integration constant. The integration constant is subsumed thanks to a parameter change $s \to s - c$. Then Eq. (8) is written as $\pm \tan\theta \frac{r}{\sqrt{1+r^2}} = -s$, which gives us

$$r(s) = \pm \frac{s}{\sqrt{\tan^2\theta - s^2}}.$$

Conversely, assume that $r(s) = \pm \dfrac{s}{\sqrt{\tan^2\theta - s^2}}$ holds and let us put

$$x = \mp \frac{\sin\theta}{\sqrt{1+r^2}} = \mp \frac{\sin\theta}{\sqrt{1+\dfrac{s^2}{\tan^2\theta - s^2}}} = \mp \cos\theta \sqrt{\tan^2\theta - s^2}, \tag{9}$$

where we are assuming that when $r$ has the positive (negative) sign, then $x$ gets the negative (positive) sign and $\theta$ is constant. Thus, $rx = -s\cos\theta$. Let now consider the vector $\vec{d}$ defined by $\vec{d} = \cos\theta\left(s\vec{t} + \vec{n}_1 \pm \left(\sqrt{\tan^2\theta - s^2}\right)\vec{n}_2\right)$. By differentiating last equality and using Frenet formulae we have $\vec{d}' = 0$, i.e., the direction of $\vec{d}$ is constant and $h(\vec{d}, \vec{n}_1) = \cos\theta = const$.

Once the intrinsic or natural equations of a curve have been determined, the next step is to integrate Frenet's formulae with $k \equiv 1$ and

$$r(s) = \pm \frac{s}{\sqrt{\tan^2\theta - s^2}} = \pm \frac{\dfrac{s}{\tan\theta}}{\sqrt{1 - \left(\dfrac{s}{\tan\theta}\right)^2}} = \pm \tan\left(\arcsin\left(\frac{s}{\tan\theta}\right)\right). \tag{10}$$

This completes the proof.

***Theorem 4.*** *The spatial quaternionic curves with principal curvature $k \equiv 1$ and such that their normal vectors make a constant angle with a fixed line are, up to rigid movements in space or up to the antipodal map, Salkowski curves.*

**Proof:** From Definition 3, the arc-length parameter of Salkowski curve is $s = \dfrac{\sin(nt)}{m}$. Therefore, $t = \dfrac{\arcsin(ms)}{n}$. Then curvature and torsion of Salkowski curve are $k \equiv 1$, $r(s) = \tan(\arcsin(ms))$, respectively. The same intrinsic equations, with $m = 1/\tan\theta$ and $n = \cos\theta$ as the ones shown in Theorem 3. For the negative case in Eq. (10), let us recall that if a curve $\alpha$ has torsion $r^\alpha$ then the curve $\beta(t) = -\alpha(t)$ has as torsion $r^\beta(t) = -r^\alpha(t)$, whereas curvature is preserved. Therefore, the fundamental theorem of spatial quaternionic curves states that up to rigid movements or up to antipodal map $p \to -p$, the curves are spatial quaternionic Salkowski curves. □

From Theorem 3 and Theorem 4, we have the following corollary.

***Corollary 1.*** *Let $\alpha(s): I \subset IR \to Q$ be an arc-length parametrized spatial quaternionic curve in $E^3$ with principal curvature $k \equiv 1$. Then $\alpha(s)$ is a spatial quaternionic slant helix if and only if $\alpha(s)$ is a spatial quaternionic Salkowski curve.*

As a special case, from Theorem 3, we have the following corollary.

***Corollary 2.*** *Let $\alpha(s): I \subset IR \to Q$ be an arc-length parametrized spatial quaternionic curve in $E^3$ with principal curvature $k \equiv 1$. The principal normal vector $\vec{n}_1(s)$ makes a constant angle $\theta = \pm\arccos(a)$, $a = const.$ with a fixed direction if and only if $r(s) = \pm \dfrac{s}{\sqrt{b^2 - s^2}}$ holds, where $b = \dfrac{a}{\sqrt{1-a^2}}$.*

## 3. Spatial Quaternionic Anti-Salkowski Curves in $E^3$

In this section, we build a spatial quaternionic curve with constant torsion from a spatial quaternionic curve of constant curvature.

Let us recall that the curve $\alpha: I \subset IR \to Q$ be a spatial quaternionic curve in $E^3$ with $\alpha' \neq 0$, $k \neq 0$. Then we have the followings:

**Theorem 5.** *Let $\alpha(s): I \subset IR \to Q$ be an arc-length parametrized spatial quaternionic curve in $E^3$ with curvatures $k^\alpha$, $r^\alpha$ and quaternionic frame $\{\vec{t}^\alpha, \vec{n}_1^\alpha, \vec{n}_2^\alpha\}$. Let us consider the spatial quaternionic curve*

$$\beta(s) = \int_{s_0}^{s} \vec{n}_2^\alpha(u)du, \tag{11}$$

*with curvatures $k^\beta$, $r^\beta$ and quaternionic frame $\{\vec{t}^\beta, \vec{n}_1^\beta, \vec{n}_2^\beta\}$. Then there exist the following relationships between the curvatures and frames of $\alpha$ and $\beta$,*

$$k^\beta = |r^\alpha|, \quad r^\beta = k^\alpha, \quad \vec{t}^\beta = \vec{n}_2^\alpha, \quad \vec{n}_1^\beta = \vec{n}_1^\alpha, \quad \vec{n}_2^\beta = -\vec{t}^\alpha.$$

**Proof:** From Eq. (11) we have $\dfrac{d\vec{\beta}}{ds} = \vec{n}_2^\alpha$. Since $\vec{n}_2^\alpha$ is unit, we get $h\left(\dfrac{d\vec{\beta}}{ds}, \dfrac{d\vec{\beta}}{ds}\right) = 1$, i.e., $\beta$ is a unit speed spatial quaternionic curve with arc-length $s$ and $\vec{t}^\beta(s) = \dfrac{d\vec{\beta}}{ds} = \vec{n}_2^\alpha(s)$. Differentiating that with respect to $s$, it follows

$$\frac{\vec{t}^\beta}{ds} = \frac{d\vec{n}_2^\alpha}{ds} = -r^\alpha \vec{n}_1^\alpha. \tag{12}$$

Therefore, $k^\beta(s) = \left\|\dfrac{d\vec{t}^\beta}{ds}\right\| = |r^\alpha|$. Then from the Frenet formulae and Eq. (12), it follows $\vec{n}_1^\beta = \vec{n}_1^\alpha$. Finally, we get $\vec{n}_2^\beta = \vec{t}^\beta \times \vec{n}_1^\beta = \vec{n}_2^\alpha \times \vec{n}_1^\alpha = -\vec{t}^\alpha$. □

Let now consider Theorem 5 for the spatial quaternionic Salkowski curve $\alpha_m(t)$ given above. From Eq. (11) we write $\vec{\beta} = \int \vec{n}_2^\alpha(s)ds = \int \vec{n}_2^\alpha(t)\|\vec{\alpha}_m'(t)\|dt$. Thus the parametric equation of the quaternionic curve $\beta$ is

$$\begin{aligned}
\beta_m(t) &= \frac{1}{2(4n^2-1)m}\left(n(1-4n^2+3\cos(2nt))\cos t + (2n^2+1)\sin t \sin(2nt)\right)\vec{e}_1 \\
&+ \frac{1}{2(4n^2-1)m}\left(n(1-4n^2+3\cos(2nt))\sin t - (2n^2+1)\cos t \sin(2nt)\right)\vec{e}_2 \\
&+ \frac{n^2-1}{4n^2}\left(2nt + \sin(2nt)\right)\vec{e}_3
\end{aligned} \tag{13}$$

where $n = \dfrac{m}{\sqrt{1+m^2}}$. Then the curvature and torsion of $\beta_m(t)$ are $k^\beta = |\tan(nt)|$, $r^\beta = 1$, respectively. After these computations we give the followings.

**Definition 4.** For any $m \in IR$, the spatial quaternionic curve $\beta_m(t)$ given in Eq. (13) with constant torsion and non-constant curvature is called spatial quaternionic anti-Salkowski curve in $E^3$.

**Theorem 6.** *Spatial Quaternionic anti-Salkowski curve $\beta_m(t)$ given in Eq. (13) is a quaternionic slant helix.*

**Proof:** Let $\beta_m(t)$ be a spatial quaternionic anti-Salkowski curve with constant torsion $r^\beta = 1$. From Theorem 5, we know that $\vec{n}_1^\beta = \vec{n}_1^\alpha$. Then from Theorem 4, $\beta_m(t)$ is a spatial quaternionic slant helix. □

## 4. Quaternionic Similar Curves with Variable Transformation in $E^3$

In this section, we give the definition and characterizations of quaternionic similar curves with variable transformation. Before giving the characterizations, first we give following definition and theorem.

**Definition 5.** Let $\alpha(s_\alpha)$ and $\beta(s_\beta)$ be two spatial quaternionic curves in $E^3$ parameterized by arc-lengths $s_\alpha$, $s_\beta$ and let curvatures, torsions and Frenet frames of the curves be $k_\alpha$, $k_\beta$; $r_\alpha$, $r_\beta$ and $\{\vec{t}^\alpha, \vec{n}_1^\alpha, \vec{n}_2^\alpha\}$, $\{\vec{t}^\beta, \vec{n}_1^\beta, \vec{n}_2^\beta\}$, respectively. Then, $\alpha(s_\alpha)$ and $\beta(s_\beta)$ are called quaternionic similar curves with variable transformation $\lambda_\beta^\alpha$ if there exists a variable transformation $s_\alpha = \int \lambda_\beta^\alpha(s_\beta) ds_\beta$ of the arc-lengths such that the tangent vectors are the same for two curves i.e., $\vec{t}^\alpha = \vec{t}^\beta$ for all corresponding values of parameters under the transformation $\lambda_\beta^\alpha$. All curves satisfying this condition is called a family of quaternionic similar curves.

**Theorem 7.** *Let $\alpha(s)$ be a spatial quaternionic curve in $E^3$ parameterized by arc-length parameter $s$. Suppose that $\alpha(\varphi)$ be another parametrization of the curve with parameter $\varphi = \int k(s) ds$. Then unit tangent vector $\vec{t}$ of $\alpha(s)$ satisfies a vector differential equation of third order given by*

$$\frac{d}{d\varphi}\left(\frac{1}{f(\varphi)}\frac{d^2\vec{t}}{d\varphi^2}\right) + \left(\frac{1+f^2(\varphi)}{f(\varphi)}\right)\frac{d\vec{t}}{d\varphi} - \left(\frac{1}{f^2(\varphi)}\frac{df(\varphi)}{d\varphi}\right)\vec{t} = 0, \qquad (14)$$

*where $f(\varphi) = \dfrac{r(\varphi)}{k(\varphi)}$.*

**Proof:** If we write derivatives given in Eq. (5) according to $\varphi$, we have

$$\frac{d\vec{t}}{d\varphi} = \vec{n}_1, \quad \frac{d\vec{n}_1}{d\varphi} = -\vec{t} + f(\varphi)\vec{n}_2, \quad \frac{d\vec{n}_2}{d\varphi} = -f(\varphi)\vec{n}_1, \qquad (15)$$

respectively, where $f(\varphi) = \dfrac{r(\varphi)}{k(\varphi)}$. From the first and second equations of (15), we get $\vec{n}_2 = \dfrac{1}{f(\varphi)}\left(\dfrac{d^2\vec{t}}{d\varphi^2} + \vec{t}\right)$ and substituting that in the last equation of (15) we have desired equation (14). □

Now we give the following theorems characterizing spatial quaternionic similar curves. In the following theorems, whenever we talk about $\alpha(s_\alpha)$ and $\beta(s_\beta)$ we assume that these curves are considered as given in Definition 5.

**Theorem 8.** Let $\alpha(s_\alpha)$ and $\beta(s_\beta)$ be two spatial quaternionic curves in $E^3$. Then $\alpha(s_\alpha)$ and $\beta(s_\beta)$ are quaternionic similar curves with variable transformation if and only if the principal normal vectors of the curves are the same, i.e.,

$$\vec{n}_1^\alpha(s_\alpha) = \vec{n}_1^\beta(s_\beta), \tag{16}$$

under the particular variable transformation

$$\lambda_\beta^\alpha = \frac{ds_\alpha}{ds_\beta} = \frac{k_\beta}{k_\alpha}, \tag{17}$$

of the arc-lengths.

**Proof:** Let $\alpha(s_\alpha)$ and $\beta(s_\beta)$ be spatial quaternionic similar curves with variable transformation. Then differentiating the equality $\vec{t}^\alpha = \vec{t}^\beta$ with respect to $s_\beta$ it follows $k_\alpha \lambda_\beta^\alpha \vec{n}_1^\alpha = k_\beta \vec{n}_1^\beta$, and we obtain Eq. (16) and (17).

Conversely, let $\alpha(s_\alpha)$ and $\beta(s_\beta)$ be two spatial quaternionic curves satisfying Eq. (16) and (17). By multiplying Eq. (16) with $k_\beta$ and differentiating the result equality with respect to $s_\beta$ we have

$$\int k_\beta(s_\beta)\vec{n}_1^\beta(s_\beta)ds_\beta = \int k_\beta(s_\beta)\vec{n}_1^\beta(s_\beta)\frac{ds_\beta}{ds_\alpha}ds_\alpha. \tag{18}$$

From Eq. (16), (17) and (18) we obtain

$$\vec{t}^\beta(s_\beta) = \int k_\beta(s_\beta)\vec{n}_1^\beta(s_\beta)ds_\beta = \int k_\alpha(s_\alpha)\vec{n}_1^\alpha(s_\alpha)ds_\alpha = \vec{t}^\alpha(s_\alpha), \tag{19}$$

which means that $\alpha(s_\alpha)$ and $\beta(s_\beta)$ are spatial quaternionic similar curves with variable transformation. □

**Theorem 9.** Let $\alpha(s_\alpha)$ and $\beta(s_\beta)$ be two spatial quaternionic curves. Then $\alpha(s_\alpha)$ and $\beta(s_\beta)$ are spatial quaternionic similar curves with variable transformation if and only if the binormal vectors of the curves are the same, i.e.,

$$\vec{n}_2^\alpha(s_\alpha) = \vec{n}_2^\beta(s_\beta), \tag{20}$$

under the particular variable transformation

$$\lambda_\beta^\alpha = \frac{ds_\alpha}{ds_\beta} = \frac{r_\beta}{r_\alpha}, \tag{21}$$

of the arc-lengths.

**Proof:** From Definition 5 and Theorem 8, there exists a variable transformation of the arc-lengths such that the tangent vectors and principal normal vectors are the same. Then we have $\vec{n}_2^\alpha(s_\alpha) = \vec{t}^\alpha(s_\alpha) \times \vec{n}_1^\alpha(s_\alpha) = \vec{t}^\beta(s_\beta) \times \vec{n}_1^\beta(s_\beta) = \vec{n}_2^\beta(s_\beta)$.

Conversely, let $\alpha(s_\alpha)$ and $\beta(s_\beta)$ be two spatial quaternionic curves satisfying Eq. (20) and (21). By differentiating Eq. (20) with respect to $s_\beta$ we obtain $r_\alpha(s_\alpha)\vec{n}_1^\alpha(s_\alpha)\frac{ds_\alpha}{ds_\beta} = r^\beta(s_\beta)\vec{n}_1^\beta(s_\beta)$ which gives us

$$\lambda_\beta^\alpha = \frac{r_\beta}{r_\alpha}, \quad \vec{n}_1^\alpha(s_\alpha) = \vec{n}_1^\beta(s_\beta). \tag{22}$$

Then from Eq. (20) and (22), we have $\vec{t}^\alpha(s_\alpha) = \vec{n}_1^\alpha(s_\alpha) \times \vec{n}_2^\alpha(s_\alpha) = \vec{n}_1^\beta(s_\beta) \times \vec{n}_2^\beta(s_\beta) = \vec{t}^\beta(s_\beta)$. □

**Theorem 10.** Let $\alpha(s_\alpha)$ and $\beta(s_\beta)$ be two spatial quaternionic curves. Then $\alpha(s_\alpha)$ and $\beta(s_\beta)$ are spatial quaternionic similar curves with variable transformation if and only if the ratio of curvatures are the same, i.e.,

$$\frac{r_\beta(s_\beta)}{k_\beta(s_\beta)} = \frac{r_\alpha(s_\alpha)}{k_\alpha(s_\alpha)}, \tag{23}$$

under the particular variable transformation keeping equal total curvatures, i.e.,

$$\varphi_\beta(s_\beta) = \int k_\beta(s_\beta)ds_\beta = \int k_\alpha(s_\alpha)ds_\alpha = \varphi_\alpha(s_\alpha) \tag{24}$$

of the arc-lengths.

**Proof:** Let $\alpha(s_\alpha)$ and $\beta(s_\beta)$ be two spatial quaternionic similar curves with variable transformation. Then from Eq. (21) and (22), we have Eq. (23) under the variable transformation Eq. (24), and this transformation is also obtained by integrating Eq. (21).

Conversely, let $\alpha(s_\alpha)$ and $\beta(s_\beta)$ be two spatial quaternionic curves satisfying Eq. (23) and (24). From Theorem 7, the unit tangents $\vec{t}^\alpha$ and $\vec{t}^\beta$ of the curves satisfy the following vector differential equations of third order

$$\frac{d}{d\varphi_\alpha}\left(\frac{1}{f_\alpha(\varphi_\alpha)}\frac{d^2\vec{t}^\alpha}{d\varphi_\alpha^2}\right) + \left(\frac{1+f_\alpha^2(\varphi_\alpha)}{f_\alpha(\varphi_\alpha)}\right)\frac{d\vec{t}^\alpha}{d\varphi_\alpha} - \left(\frac{1}{f_\alpha^2(\varphi_\alpha)}\frac{df_\alpha(\varphi_\alpha)}{d\varphi_\alpha}\right)\vec{t}^\alpha = 0, \tag{25}$$

$$\frac{d}{d\varphi_\beta}\left(\frac{1}{f_\beta(\varphi_\beta)}\frac{d^2\vec{t}^\beta}{d\varphi_\beta^2}\right) + \left(\frac{1+f_\beta^2(\varphi_\beta)}{f_\beta(\varphi_\beta)}\right)\frac{d\vec{t}^\beta}{d\varphi_\beta} - \left(\frac{1}{f_\beta^2(\varphi_\beta)}\frac{df_\beta(\varphi_\beta)}{d\varphi_\beta}\right)\vec{t}^\beta = 0, \tag{26}$$

where $f_\alpha(\varphi_\alpha) = \frac{r_\alpha(\varphi_\alpha)}{k_\alpha(\varphi_\alpha)}$, $f_\beta(\varphi_\beta) = \frac{r_\beta(\varphi_\beta)}{k_\beta(\varphi_\beta)}$, $\varphi_\alpha(s_\alpha) = \int k_\alpha(s_\alpha)ds_\alpha$, $\varphi_\beta(s_\beta) = \int k_\beta(s_\beta)ds_\beta$.

From Eq. (23), we have $f_\alpha(\varphi_\alpha) = f_\beta(\varphi_\beta)$ under the variable transformation $\varphi_\alpha = \varphi_\beta$. Thus under Eq. (23) and transformation Eq. (24), Eq. (25) and (26) are the same, i.e., they have the same solutions. It means that the unit tangents $\vec{t}^\alpha$ and $\vec{t}^\beta$ are the same. Then $\alpha(s_\alpha)$ and $\beta(s_\beta)$ are two spatial quaternionic similar curves with variable transformation. □

Let now consider some special cases. From Eq. (17) and (22), we have $k_\beta = \lambda_\beta^\alpha k_\alpha$, $r_\beta = \lambda_\beta^\alpha r_\alpha$, respectively. From last equalities it is clear that if $\alpha(s_\alpha)$ is a straight line i.e., $k_\alpha = 0$, then under the variable transformation the curvature does not change. So we have the following corollaries.

***Corollary 3.*** *The family of spatial quaternionic straight lines forms a family of spatial quaternionic similar curves with variable transformation.*

If $\alpha(s_\alpha)$ is a spatial quaternionic plane curve i.e., $r_\alpha = 0$, then under the variable transformation the curvature does not change. So we have the following corollary.

***Corollary 4.*** *The family of spatial quaternionic plane curves forms a family of spatial quaternionic similar curves with variable transformation.*

If $\alpha(s_\alpha)$ is a spatial quaternionic $\vec{n}_2^\alpha$-slant helix, then $\frac{r_\alpha}{k_\alpha} = \tan\theta$ is constant where $\theta$ is the constant angle between $\vec{n}_2^\alpha$ and axis of $\vec{n}_2^\alpha$-slant helix [12]. Then if the spatial quaternionic $\vec{n}_2^\alpha$-slant helix $\alpha(s_\alpha)$ and spatial quaternionic curve $\beta(s_\beta)$ are two spatial quaternionic similar curves with variable transformation, then from Theorem 10, we have $\frac{r_\beta}{k_\beta} = \tan\theta$ is constant which gives following corollary.

**Corollary 5.** *The family of spatial quaternionic $\vec{n}_2^\alpha$-slant helices with fixed constant angle $\theta$ between $\vec{n}_2^\alpha$ and axis of helix forms a family of spatial quaternionic similar curves with variable transformation.*

If $\alpha(s_\alpha)$ is a spatial quaternionic Salkowski curve, then from Corollary 2, we have $k_\alpha(s_\alpha) = 1$, $r_\alpha(s_\alpha) = \pm\frac{s_\alpha}{\sqrt{b^2 - s_\alpha^2}}$. Let now $\alpha(s_\alpha)$ and $\beta(s_\beta)$ be two spatial quaternionic Salkowski curves such that the transformation (24) holds. Then it is easy to prove that

$$\frac{r_\alpha(s_\alpha)}{k_\alpha(s_\alpha)} = \pm\frac{s_\alpha}{\sqrt{b^2 - s_\alpha^2}} = \pm\frac{s_\beta}{\sqrt{b^2 - s_\beta^2}} = \frac{r_\beta(s_\beta)}{k_\beta(s_\beta)}, \quad (27)$$

which leads to the following corollary.

**Corollary 6.** *The family of quaternionic Salkowski curves with $k_\alpha(s_\alpha) = 1$ in $E^3$ forms a family of spatial quaternionic similar curves with variable transformation.*

If $\alpha(s_\alpha)$ is a spatial quaternionic anti-Salkowski curve in $E^3$, then from Corollary 2 and Theorem 5, we have $k_\alpha(s_\alpha) = \frac{s_\alpha}{\sqrt{b^2 - s_\alpha^2}}$, $r_\alpha(s_\alpha) = 1$. Let now $\alpha(s_\alpha)$ and $\beta(s_\beta)$ be two spatial quaternionic anti-Salkowski curves such that the transformation (24) holds. Then one can easily prove that

$$\frac{r_\alpha(s_\alpha)}{k_\alpha(s_\alpha)} = \frac{\sqrt{b^2 - s_\alpha^2}}{s_\alpha} = \frac{\sqrt{b^2 - s_\beta^2}}{s_\beta} = \frac{r_\beta(s_\beta)}{k_\beta(s_\beta)}, \quad (28)$$

which leads to the following corollary.

**Corollary 7.** *The family of spatial quaternionic anti-Salkowski curves with $r_\alpha(s_\alpha) = 1$ in $E^3$ forms a family of spatial quaternionic similar curves with variable transformation.*

### 5. Quaternionic Salkowski Curves in $E^4$

In this section we give the definition and characterizations of quaternionic Salkowski curves in $E^4$.

**Definition 6.** Let $\alpha : I \subset IR \to Q$ be a real spatial quaternionic curve in $E^3$ such that $\alpha(s) = \sum_{i=1}^{3} \alpha_i(s) e_i$ and let $\gamma : I \subset IR \to Q$, $s \to \gamma(s) = \sum_{i=1}^{4} \alpha_i(s) e_i$, $e_4 = +1$, be an arc-length parametrized quaternionic curve in $E^4$ obtained from $\alpha(s)$ with nonzero curvatures $\{K, k, r - K\}$. Let $\{\vec{T}(s), \vec{N}(s), \vec{B}_1(s), \vec{B}_2(s)\}$ denotes the quaternionic Frenet frame along the curve $\gamma(s)$. We call $\gamma(s)$ as quaternionic Salkowski curve in $E^4$ if the principal curvature is $K \equiv 1$ and torsion $k(s)$ and bitorsion $(r - K)(s)$ are non-constants.

**Theorem 11.** *Let $\gamma: I \subset \mathbb{R} \to Q$ be a unit speed quaternionic Salkowski curve in $E^4$. Then $\gamma$ is a quaternionic $B_2$-slant helix if and only if $\frac{r-1}{k} = \pm \tan\theta \sin s$ holds where $\theta$ is the constant angle between the vector $\vec{B}_2$ and a fixed quaternion.*

**Proof.** Let $\gamma$ be a unit speed quaternionic Salkowski curve in $E^4$. Assume that $\gamma$ is a $B_2$-slant helix. Then the second binormal unit vector $\vec{B}_2$ makes a constant angle $\theta$ with a fixed direction in a unit quaternion $\vec{U}$; that is

$$h(\vec{B}_2, \vec{U}) = \cos\theta = constant, \tag{29}$$

along the curve. By differentiation Eq. (29) with respect to $s$ and using the Frenet formulae (5) we have $(r-1)h(\vec{B}_1, \vec{U}) = 0$. Therefore $\vec{U}$ is in the subspace $Sp\{\vec{T}, \vec{N}, \vec{B}_2\}$. Then we have

$$\vec{U} = a_1 \vec{T} + a_2 \vec{N} + a_3 \vec{B}_2, \tag{30}$$

where $a_1 = a_1(s) = \cos\theta_1(s)$, $a_2 = a_2(s) = \cos\theta_2(s)$ and $a_3 = a_3(s) = \cos\theta(s)$ are direction cosines of $\vec{U}$ and from Eq. (29) $a_3 = \cos\theta$ is constant. Since $\vec{U}$ is unit, we have

$$a_1^2 + a_2^2 + a_3^2 = 1. \tag{31}$$

The differentiation of Eq. (30) gives

$$\left(\frac{da_1}{ds} - a_2\right)\vec{T} + \left(\frac{da_2}{ds} + a_1\right)\vec{N} + (a_2 k - a_3(r-1))\vec{B}_1 + \frac{da_3}{ds}\vec{B}_2 = 0,$$

and from this equation we get

$$a_2 = \frac{r-1}{k} a_3 = \frac{da_1}{ds}, \quad \frac{da_2}{ds} = -a_1. \tag{32}$$

Since $\frac{da_2}{ds} = -a_1$ and $\frac{da_2}{ds} = \frac{d^2 a_1}{ds^2}$, we find the second order linear differential equation in $a_1$ given by

$$\frac{d^2 a_1}{ds^2} + a_1 = 0, \tag{33}$$

which has the solution $a_1 = A\cos(s) + B\sin(s)$ where $A$ and $B$ are constants. From Eq. (32) and (33) we have

$$a_2 = \frac{r-1}{k} a_3 = -A\sin(s) + B\cos(s), \quad a_1 = -\left(\frac{r-1}{k}\right)' a_3 = A\cos(s) + B\sin(s).$$

From these equations it follows that

$$A = -a_3\left(\frac{r-1}{k}\sin(s) + \left(\frac{r-1}{k}\right)'\cos(s)\right), \quad B = a_3\left(\frac{r-1}{k}\cos(s) - \left(\frac{r-1}{k}\right)'\sin(s)\right), \tag{34}$$

Hence using Eq. (31) and (34) we get

$$A^2 + B^2 = \left[\left(\frac{r-1}{k}\right)^2 + \left(\left(\frac{r-1}{k}\right)'\right)^2\right]\cos^2\theta = \sin^2\theta,$$

or

$$\left(\frac{r-1}{k}\right)' = \pm\sqrt{\tan^2\theta - \left(\frac{r-1}{k}\right)^2}. \tag{35}$$

Integrating Eq. (35) gives us $\frac{r-1}{k} = \tan\theta \sin(\pm(s+c))$. Making the parameter change $s \to s-c$, the last equality is written as

$$\frac{r-1}{k} = \pm\tan\theta \sin s. \tag{36}$$

Conversely, if the condition (36) is satisfied for a quaternionic curve $\gamma$ we always find a constant vector $\vec{U}$ which makes a constant angle with the second binormal of the curve. Consider the unit vector $\vec{U}$ defined by $\vec{U} = \left[ -\left(\frac{r-1}{k}\right)' \vec{T} + \frac{r-1}{k}\vec{N} + \vec{B}_2 \right]\cos\theta$ where $\theta$ is constant. Differentiating $\vec{U}$ and using Eq. (36) gives that $\vec{U}' = 0$, this means that $\vec{U}$ is a constant quaternion. Moreover, $h(\vec{U}, \vec{B}_2) = \cos\theta$ is constant, i.e, $\gamma$ is a quaternionic $\vec{B}_2$-slant helix.

Let now consider that $\gamma(s)$ is a quaternionic Salkowski curve in $E^4$ obtained from 3-dimensional Salkowski curve $\alpha(s)$. Then considering Theorem 3, we have the following corollary.

**Corollary 8.** *Let $\gamma(s)$ be a quaternionic Salkowski curve in $E^4$ obtained from 3-dimensional Salkowski curve $\alpha(s)$. Then $\gamma$ is a quaternionic $B_2$-slant helix if and only if*

$$\pm \frac{s_\alpha}{\sqrt{\tan^2\theta_\alpha - s_\alpha^2}} = 1 \pm \tan\theta_\gamma \sin s_\gamma,$$

*holds where $s_\alpha$, $s_\gamma$ are arc-length parameters of the curves $\alpha(s)$ and $\gamma(s)$, respectively and $\theta_\alpha$ is the constant angle between the principal normal vector $\vec{n}_1(s)$ of $\alpha(s)$ and a fixed direction, $\theta_\gamma$ is the constant angle between the vector $\vec{B}_2$ of $\gamma(s)$ and a fixed direction.*

## 6. Quaternionic Anti-Salkowski Curves in $E^4$

In this section we build a quaternionic curve in $E^4$ with constant bitorsion from a quaternionic curve of constant principal curvature.

Let us recall that the curve $\gamma: I \subset IR \to Q$ be a quaternionic curve in $E^4$ with $\gamma' \neq 0$, $K \neq 0$. Then we have the following theorem:

**Theorem 12.** *Let $\gamma(s): I \subset IR \to Q$ be an arc-length parametrized quaternionic curve in $E^4$ with curvatures $K^\gamma$, $k^\gamma$, $(r-K)^\gamma$ and quaternionic frame $\{\vec{T}^\gamma(s), \vec{N}^\gamma(s), \vec{B}_1^\gamma(s), \vec{B}_2^\gamma(s)\}$. Let us consider the quaternionic curve*

$$\beta(s) = \int_{s_0}^{s} \vec{B}_2^\gamma(u)\,du, \tag{37}$$

*in $E^4$ with curvatures $K^\beta$, $k^\beta$, $(r-K)^\beta$ and quaternionic frame $\{\vec{T}^\beta, \vec{N}^\beta, \vec{B}_1^\beta, \vec{B}_2^\beta\}$. Then there exist the following relationships between the curvatures and Frenet vectors of $\gamma$ and $\beta$,*

$$K^\beta = \left|(r-K)^\gamma\right|,\ k^\beta = \left|k^\gamma\right|,\ K^\gamma = \left|(r-K)^\beta\right|,\ \vec{T}^\beta = \vec{B}_2^\gamma,\ \vec{N}^\beta = \vec{B}_1^\gamma,\ \vec{B}_1^\beta = \vec{N}^\gamma,\ \vec{B}_2^\beta = \vec{T}^\gamma. \tag{38}$$

**Proof:** From Eq. (37) we have $\frac{d\vec{\beta}}{ds} = \vec{B}_2^\gamma$. Since $\vec{B}_2^\gamma$ is unit, we get $h\left(\frac{d\vec{\beta}}{ds}, \frac{d\vec{\beta}}{ds}\right) = 1$, i.e., $\beta$ is a unit speed quaternionic curve with arc length $s$ and

$$\vec{T}^\beta(s) = \frac{d\vec{\beta}}{ds} = \vec{B}_2^\gamma(s). \tag{39}$$

Differentiating Eq. (39) with respect to $s$, it follows

$$\frac{d\vec{T}^\beta}{ds} = \frac{d\vec{B}_2^\gamma}{ds} = -(r-K)^\gamma \vec{B}_1^\gamma. \tag{40}$$

Therefore,

$$K^\beta(s) = \left\|\frac{d\vec{T}^\beta}{ds}\right\| = \left|(r-K)^\gamma\right|. \tag{41}$$

Then from the Frenet formulae and Eq. (40), we have

$$\vec{N}^\beta = \vec{B}_1^\gamma. \tag{42}$$

Differentiating Eq. (42) gives $-K^\beta \vec{T}^\beta + k^\beta \vec{B}_1^\beta = -k^\gamma \vec{N}^\gamma + (r-K)^\gamma \vec{B}_2^\gamma$ and using Eq. (39) and (41) we have $k^\beta = |k^\gamma|$, $K^\gamma = \left|(r-K)^\beta\right|$ and $\vec{B}_1^\beta = \vec{N}^\gamma$. □

Let now consider Theorem 12 for the quaternionic Salkowski curve $\gamma(s)$ defined in Definition 6. Then we have $(r-K)^\beta = \pm 1$, and we give the following definition.

**Definition 7.** The quaternionic curve $\beta(s)$ with constant bitorsion $(r-K)^\beta = \pm 1$, non-constant principal curvature $K^\beta$ and non-constant torsion $k^\beta$ is called quaternionic anti-Salkowski curve.

Then we give the following corollary whose proof is obtained from Theorem 12.

***Corollary 9.*** *Let quaternionic curve $\beta(s)$ be an anti-Salkowski curve in $E^4$ obtained from quaternionic Salkowski curve $\gamma(s)$. Then*
  *(i) $\beta(s)$ is a helix if and only if $\gamma(s)$ is a $B_2$-slant helix.*
  *(ii) $\beta(s)$ is a slant helix if and only if $\gamma(s)$ is a $B_1$-slant helix.*
  *(iii) $\beta(s)$ is a $B_1$-slant helix if and only if $\gamma(s)$ is a slant helix.*
  *(iv) $\beta(s)$ is a $B_2$-slant helix if and only if $\gamma(s)$ is a helix.*

## 7. Quaternionic Similar Curves with Variable Transformation in $E^4$

In this section we give the definition and characterizations of quaternionic similar curves with variable transformation in $E^4$. Before giving the characterizations, first we give the following definition and theorem.

**Definition 8.** Let $\alpha(s_\alpha)$ and $\beta(s_\beta)$ be two quaternionic curves in $E^4$ parameterized by arc-lengths $s_\alpha$, $s_\beta$ with non-zero curvatures $K_\alpha$, $k_\alpha$, $r_\alpha - K_\alpha$; $K_\beta$, $k_\beta$, $r_\beta - K_\beta$ and Frenet frames $\{\vec{T}^\alpha(s_\alpha), \vec{N}^\alpha(s_\alpha), \vec{B}_1^\alpha(s_\alpha), \vec{B}_2^\alpha(s_\alpha)\}$ and $\{\vec{T}^\beta(s_\beta), \vec{N}^\beta(s_\beta), \vec{B}_1^\beta(s_\beta), \vec{B}_2^\beta(s_\beta)\}$, respectively. Then, $\alpha(s_\alpha)$ and $\beta(s_\beta)$ are called quaternionic similar curves with variable transformation $\lambda_\beta^\alpha$ if there exists a variable transformation $s_\alpha = \int \lambda_\beta^\alpha(s_\beta) ds_\beta$ of the arc-lengths such that the tangent quaternions are the same for two curves i.e., $\vec{T}^\alpha = \vec{T}^\beta$, for all corresponding values of parameters under the transformation $\lambda_\beta^\alpha$. The set of all quaternionic curves satisfying this condition is called a family of quaternionic similar curves in $E^4$.

If we integrate the equality $\vec{T}^\alpha = \vec{T}^\beta$, we have the following theorem.

**Theorem 13.** *The position vectors of the family of quaternionic similar curves with variable transformation can be written in the following form,*

$$\beta(s_\beta) = \int \vec{T}^\alpha \left(s_\alpha(s_\beta)\right) ds_\beta = \int \vec{T}^\beta(s_\alpha) \lambda_\alpha^\beta ds_\beta.$$

Before giving the characterizations of quaternionic similar curves in $E^4$, we first give the following theorem.

**Theorem 14.** *Let $\alpha(s)$ be a quaternionic curve in $E^4$ parameterized by arc-length $s$. Suppose that $\alpha(\varphi)$ be another parametrization of the curve with parameter $\varphi = \int K(s) ds$. Then the unit tangent vector $\vec{T}$ of $\alpha(s)$ satisfies a vector differential equation of fourth order given by*

$$\left[\frac{1}{g}\left(\left[\frac{1}{f}(\vec{T}''+\vec{T})\right]' + f\vec{T}'\right)\right]' + \frac{g}{f}(\vec{T}''+\vec{T}) = 0, \tag{43}$$

where $f(\varphi) = \dfrac{k(\varphi)}{K(\varphi)}$, $g(\varphi) = \dfrac{r(\varphi)}{K(\varphi)} - 1$ and $\vec{T}' = \dfrac{d\vec{T}}{d\varphi}$, $\vec{T}'' = \dfrac{d^2\vec{T}}{d\varphi^2}$.

**Proof:** If we write derivatives given in Eq. (5) according to $\varphi$, we have

$$\begin{cases} \dfrac{d\vec{T}}{d\varphi} = (K\vec{N}_1)\dfrac{1}{K} = \vec{N}_1, \\[4pt] \dfrac{d\vec{N}}{d\varphi} = (-K\vec{T} + k\vec{B}_1)\dfrac{1}{K} = -\vec{T} + f\vec{B}_1, \\[4pt] \dfrac{d\vec{B}_1}{d\varphi} = \left(-k\vec{N} + (r-K)\vec{B}_2\right)\dfrac{1}{K} = -f\vec{N} + g\vec{B}_2, \\[4pt] \dfrac{d\vec{B}_2}{d\varphi} = \left(-(r-K)\vec{B}_1\right)\dfrac{1}{K} = -g\vec{B}_1, \end{cases} \tag{44}$$

respectively, where $f(\varphi) = \dfrac{r(\varphi)}{K(\varphi)}$, $g(\varphi) = \dfrac{r(\varphi)}{K(\varphi)} - 1$. From the first and second equations of (44), we have $\vec{B}_1 = \dfrac{1}{f}(\vec{T}'' + \vec{T})$. Differentiating that with respect to $\varphi$ and using the second and third equations of (44), it follows $\vec{B}_2 = \dfrac{1}{g}\left(\left[\dfrac{1}{f}(\vec{T}'' + \vec{T})\right]' + f\vec{T}'\right)$. Finally, by differentiating last equality with respect to $\varphi$ and using the last equation of (44), we have desired equation (43) □.

Now we give the following theorems characterizing quaternionic similar curves in $E^4$. In the following theorems, whenever we talk about $\alpha(s_\alpha)$ and $\beta(s_\beta)$ we assume that these curves are defined as given in Definition 8.

**Theorem 15.** Let $\alpha(s_\alpha)$ and $\beta(s_\beta)$ be two quaternionic curves in $E^4$. Then $\alpha(s_\alpha)$ and $\beta(s_\beta)$ are quaternionic similar curves with variable transformation if and only if the second vectors of the frames are the same, i.e.,

$$\vec{N}^\alpha(s_\alpha) = \vec{N}^\beta(s_\beta), \tag{45}$$

under the particular variable transformation

$$\lambda_\beta^\alpha = \frac{K_\beta}{K_\alpha}, \tag{46}$$

of the arc-lengths.

**Proof:** Let $\alpha(s_\alpha)$ and $\beta(s_\beta)$ be quaternionic similar curves with variable transformation in $E^4$. Then differentiating the equality $\vec{T}^\alpha = \vec{T}^\beta$ with respect to $s_\beta$ it follows, $\frac{ds_\alpha}{ds_\beta} K_\alpha \vec{N}^\alpha = K_\beta \vec{N}^\beta$. Then, we obtain Eq. (45) and (46).

Conversely, let $\alpha(s_\alpha)$ and $\beta(s_\beta)$ be two quaternionic curves in $E^4$ satisfying Eq. (45) and (46). By multiplying Eq. (45) with $K_\beta$ and integrating the result equality with respect to $s_\beta$ we have

$$\int K_\beta(s_\beta)\vec{N}^\beta(s_\beta)ds_\beta = \int K_\beta(s_\beta)\vec{N}^\beta(s_\beta)\frac{ds_\beta}{ds_\alpha}ds_\alpha. \tag{47}$$

From Eq. (45), (46) and (47), we obtain

$$\vec{T}^\beta(s_\beta) = \int K_\beta(s_\beta)\vec{N}^\beta(s_\beta)ds_\beta = \int K_\alpha(s_\alpha)\vec{N}^\alpha(s_\alpha)ds_\alpha = \vec{T}^\alpha(s_\alpha), \tag{48}$$

which means that $\alpha(s_\alpha)$ and $\beta(s_\beta)$ are quaternionic similar curves with variable transformation. □

**Theorem 16.** Let $\alpha(s_\alpha)$ and $\beta(s_\beta)$ be two quaternionic curves in $E^4$. Then $\alpha(s_\alpha)$ and $\beta(s_\beta)$ are quaternionic similar curves with variable transformation if and only if the Frenet vectors $\vec{B}_1^\alpha$ and $\vec{B}_1^\beta$ of the curves are the same, i.e.,

$$\vec{B}_1^\alpha(s_\alpha) = \vec{B}_1^\beta(s_\beta), \tag{49}$$

under the particular variable transformation

$$\lambda_\beta^\alpha = \frac{K_\beta}{K_\alpha} = \frac{k_\beta}{k_\alpha}, \tag{50}$$

keeping equal total curvatures, i.e., $\varphi_\beta(s_\beta) = \int K_\beta(s_\beta)ds_\beta = \int K_\alpha(s_\alpha)ds_\alpha = \varphi_\alpha(s_\alpha)$ of the arc-lengths.

**Proof:** From Definition 8 and Theorem 14, there exists a variable transformation of the arc-lengths such that $\vec{T}^\alpha(s_\alpha) = \vec{T}^\beta(s_\beta)$ and $\vec{N}^\alpha(s_\alpha) = \vec{N}^\beta(s_\beta)$. Differentiating the second equality with respect to $s_\beta$ and considering the first one, we have

$$-K_\beta \vec{T}^\beta + k_\beta \vec{B}_1^\beta = \frac{ds_\alpha}{ds_\beta}(-K_\alpha \vec{T}^\alpha + k_\alpha \vec{B}_1^\alpha) \tag{51}$$

which gives us $\vec{B}_1^\alpha(s_\alpha) = \vec{B}_1^\beta(s_\beta)$ and $\lambda_\beta^\alpha = \frac{ds_\alpha}{ds_\beta} = \frac{K_\beta}{K_\alpha} = \frac{k_\beta}{k_\alpha}$.

Conversely, let $\alpha(s_\alpha)$ and $\beta(s_\beta)$ be two quaternionic curves in $E^4$ satisfying Eq. (49) and (50). By differentiating Eq. (49) with respect to $s_\beta$ we obtain

$$-k_\beta \vec{N}^\beta + (r_\beta - K_\beta)\vec{B}_2^\beta = \frac{ds_\alpha}{ds_\beta}\left(-k_\alpha \vec{N}^\alpha + (r_\alpha - K_\alpha)\vec{B}_2^\alpha\right). \tag{52}$$

From Eq. (50) and (52), it follows

$$-k_\beta \vec{N}^\beta + (r_\beta - K_\beta)\vec{B}_2^\beta = -k_\beta \vec{N}^\alpha + \frac{ds_\alpha}{ds_\beta}(r_\alpha - K_\alpha)\vec{B}_2^\alpha. \tag{53}$$

Since the Frenet quaternions are unit, from Eq. (53) if follows, $\dfrac{ds_\alpha}{ds_\beta} = \dfrac{r_\beta - K_\beta}{r_\alpha - K_\alpha}$. Then from Eq. (50) we have

$$\frac{r_\beta}{K_\beta} = \frac{r_\alpha}{K_\alpha} \tag{54}$$

Let now consider Theorem 14. Then, the unit tangents $\vec{T}^\alpha$ and $\vec{T}^\beta$ of the curves satisfy the following vector differential equations of fourth order

$$\left[\frac{1}{g_\alpha}\left(\left[\frac{1}{f_\alpha}\left((\vec{T}^\alpha)'' + \vec{T}^\alpha\right)\right]' + f(\vec{T}^\alpha)'\right)\right]' + \frac{g_\alpha}{f_\alpha}\left((\vec{T}^\alpha)'' + \vec{T}^\alpha\right) = 0, \tag{55}$$

$$\left[\frac{1}{g_\alpha}\left(\left[\frac{1}{f_\alpha}\left((\vec{T}^\alpha)'' + \vec{T}^\alpha\right)\right]' + f(\vec{T}^\alpha)'\right)\right]' + \frac{g_\alpha}{f_\alpha}\left((\vec{T}^\alpha)'' + \vec{T}^\alpha\right) = 0, \tag{56}$$

respectively, where

$$f_\alpha(\varphi_\alpha) = \frac{k_\alpha(\varphi_\alpha)}{K_\alpha(\varphi_\alpha)}, \quad g_\alpha(\varphi_\alpha) = \frac{r_\alpha(\varphi_\alpha)}{K_\alpha(\varphi_\alpha)} - 1, \quad f_\beta(\varphi_\beta) = \frac{k_\beta(\varphi_\beta)}{K_\beta(\varphi_\beta)},$$

$$g_\beta(\varphi_\beta) = \frac{r_\beta(\varphi_\beta)}{K_\beta(\varphi_\beta)} - 1, \quad \varphi_\alpha(s_\alpha) = \int K_\alpha(s_\alpha)ds_\alpha, \quad \varphi_\beta(s_\beta) = \int K_\beta(s_\beta)ds_\beta.$$

Now, from Eq. (50) and (54) we have $f_\alpha(\varphi_\alpha) = f_\beta(\varphi_\beta)$ and $g_\alpha(\varphi_\alpha) = g_\beta(\varphi_\beta)$ under the variable transformation $\varphi_\alpha = \varphi_\beta$. Then under the equations (50), (54) and transformation $\varphi_\alpha = \varphi_\beta$, Eq. (55) and (56) are the same, i.e., they have the same solutions. It means that the unit tangents $\vec{T}^\alpha$ and $\vec{T}^\beta$ are the same. Then $\alpha(s_\alpha)$ and $\beta(s_\beta)$ are two quaternionic similar curves with variable transformation in $E^4$. □

**Theorem 17.** Let $\alpha(s_\alpha)$ and $\beta(s_\beta)$ be two quaternionic curves in $E^4$. Then $\alpha(s_\alpha)$ and $\beta(s_\beta)$ are quaternionic similar curves with variable transformation if and only if the Frenet vectors $\vec{B}_2^\alpha$ and $\vec{B}_2^\beta$ of the curves are the same, i.e.,

$$\vec{B}_2^\alpha(s_\alpha) = \vec{B}_2^\beta(s_\beta), \tag{57}$$

*under the particular variable transformation ,*

$$\lambda_\beta^\alpha = \frac{k_\beta}{k_\alpha} = \frac{r_\beta - K_\beta}{r_\alpha - K_\alpha}, \tag{58}$$

*of the arc-lengths.*

**Proof:** Let $\alpha(s_\alpha)$ and $\beta(s_\beta)$ be two quaternionic similar curves with variable transformation. Then from Theorem 16, we have $\vec{B}_1^\alpha(s_\alpha) = \vec{B}_1^\beta(s_\beta)$. Differentiating this equality with respect to $s_\beta$ gives

$$-k_\beta \vec{N}^\beta + (r_\beta - K_\beta)\vec{B}_2^\beta = \frac{ds_\alpha}{ds_\beta}\left(-k_\alpha \vec{N}^\alpha + (r_\alpha - K_\alpha)\vec{B}_2^\alpha\right) \tag{59}$$

By considering Theorem 14, Eq. (59) gives $\frac{ds_\alpha}{ds_\beta} = \frac{k_\beta}{k_\alpha} = \frac{r_\beta - K_\beta}{r_\alpha - K_\alpha}$ and $\vec{B}_2^\alpha(s_\alpha) = \vec{B}_2^\beta(s_\beta)$.

Conversely, let $\alpha(s_\alpha)$ and $\beta(s_\beta)$ be two quaternionic curves satisfying Eq. (57) and (58). Differentiating Eq. (57) with respect to $s_\beta$ it follows $-(r_\beta - K_\beta)\vec{B}_1^\beta = -\frac{ds_\alpha}{ds_\beta}(r_\alpha - K_\alpha)\vec{B}_1^\alpha$. From Eq. (58), we see that $\vec{B}_1^\beta = \vec{B}_1^\alpha$. Then by Theorem 16, we obtain that $\alpha(s_\alpha)$ and $\beta(s_\beta)$ are two quaternionic similar curves with variable transformation.□

Let now consider some special cases. From Eq. (46) and (58) we have
$$K_\beta = \lambda_\beta^\alpha K_\alpha, \quad k_\beta = \lambda_\beta^\alpha k_\alpha, \quad r_\beta - K_\beta = \lambda_\beta^\alpha (r_\alpha - K_\alpha), \tag{60}$$
respectively. From Eq. (60) it is clear that if the principal curvature $K_\alpha$ of $\alpha(s_\alpha)$ vanishes, i.e., $K_\alpha \equiv 0$, then under the variable transformation the principal curvature $K_\beta$ of $\beta(s_\beta)$ does not change, i.e, $K_\beta \equiv 0$. Similarly, if the torsion $k_\alpha$ of $\alpha(s_\alpha)$ is given by $k_\alpha \equiv 0$, then under the variable transformation the torsion $k_\beta$ of $\beta(s_\beta)$ is also given by $k_\beta \equiv 0$. Finally, if the bitorsion $r_\alpha - K_\alpha$ of $\alpha(s_\alpha)$ vanishes, then we see that the bitorsion $r_\beta - K_\beta$ of $\beta(s_\beta)$ also vanishes. So, we have the following corollaries.

***Corollary 10.*** *The family of quaternionic curves in $E^4$ with vanishing principal curvature $K$ forms a family of quaternionic similar curves with variable transformation.*

***Corollary 11.*** *The family of quaternionic curves in $E^4$ with vanishing torsion $k$ forms a family of quaternionic similar curves with variable transformation.*

***Corollary 12.*** *The family of quaternionic curves in $E^4$ with vanishing bitorsion $r - K$ forms a family of quaternionic similar curves with variable transformation.*

## 8. Conclusions
The curves of Euclidean spaces $E^3$ and $E^4$ are considered as the spatial quaternions and quaternions and called spatial quaternionic curves and quaternionic curves, respectively. So, the types of Euclidean curves are considered for spatial quaternionic and quaternionic curves. By considering that this paper gives two new types of these curves as Salkowski curves and similar curves of the type quaternionic. The relations between these curves are obtained. Of course these curves can be studied in new spaces such as semi-Euclidean spaces and spatial quaternionic and quaternionic Salkowski curves and similar curves can be defined and studied in the spaces $E_1^3$ and $E_1^4$.